\documentclass[12pt]{amsart}
\usepackage{amsmath,amssymb,amsfonts,amscd,graphicx,xypic}

\newcommand*{\QEDB}{\hfill\ensuremath{\square}}

\newcommand{\slantedcup}{\mathbin{\rotatebox[origin=c]{50}{$\cup$}}}

\newtheorem{lemma}{Lemma}[section]

\newtheorem{proposition}[lemma]{Proposition}
\newtheorem{corollary}[lemma]{Corollary}
\theoremstyle{definition}
\newtheorem{definition}[lemma]{Definition}

\theoremstyle{remark}

\numberwithin{equation}{section}
\theoremstyle{plain}
\newtheorem{thm}{Theorem}

\theoremstyle{definition}

\newtheorem{condition}[lemma]{Condition}

\newcommand{\co}{\colon\thinspace}

\newcommand{\nl}{\hfil\break}

\numberwithin{figure}{section}

\begin{document}

\title{A homotopy$^{\mathbf +}$ solution to the A-B slice problem}
\author{Michael Freedman and Vyacheslav Krushkal}

\address{Microsoft Station Q, University of California, Santa Barbara, CA 93106-6105, and 
Department of Mathematics, University of California, Santa Barbara, CA 93106}
\email{michaelf\char 64 microsoft.com}
\address{Department of Mathematics, University of Virginia, Charlottesville, VA 22904}
\email{krushkal\char 64 virginia.edu}

\begin{abstract}
The \mbox{A-B} slice problem, a reformulation of the $4$-dimensional topological surgery conjecture for free groups, is shown to admit a link-homotopy$^+$ solution.
The proof relies on geometric applications of the group-theoretic $2$-Engel relation. 
Implications for the surgery conjecture are discussed.
\end{abstract}

\dedicatory{Dedicated to the memory of Tim Cochran}

\maketitle

\section{Introduction} \label{introduction}

Four-dimensional surgery is known to work in the topological category for a class of {\em good} fundamental groups. This result was originally established in the simply-connected case in \cite{F0}, and it is currently known to hold for groups of subexponential growth  {and a somewhat larger class generated by these} \cite{F1, FT, KQ}. The A-B slice problem \cite{F2, F3} is a reformulation of the surgery conjecture for free groups,  {which is the most difficult case}.

The A-B slice problem concerns decompositions of the $4$-ball. The handle structure of a decomposition, interpreted as a Kirby diagram, gives rise to a {\em stabilization} of a given link $L$, see section 2 for a precise definition. In these terms, to show that $L$ is A-B slice one needs to find a stabilization and band-sums between the components so that the resulting link is slice. 
 {The Generalized Borromean Rings (GBRs) are a collection of links  {any coinitial subset of which is} universal for surgery.
In a recent work \cite{FK} we showed that GBRs have a coinitial subset admitting a link-homotopy solution to the A-B slice problem.} 
In other words, given such a GBR there exists a stabilization and band sums so that the resulting link is homotopically trivial (h-trivial) in the sense of Milnor \cite{M}.
 Here we sharpen this result:

\begin{thm} \label{homotopy+ theorem} \sl
The A-B slice problem for a coinitial collection of generalized Borromean rings, forming universal surgery problems, admits a (link-homotopy)$^+$ solution.
\end{thm}

\begin{figure}[ht]
\centering
\includegraphics[height=2.7cm]{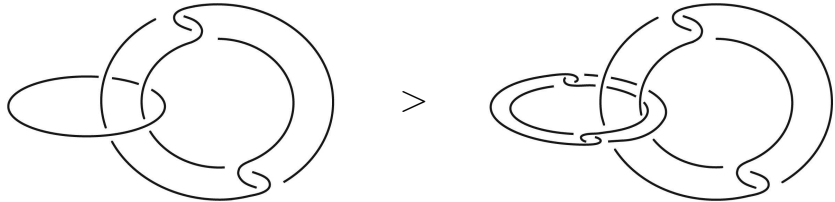}
\large{
\put(-165,35){$>$}
}
 \caption{ {Partial order on Generalized Borromean Rings}}
\label{GBRs}
\end{figure}

 {The Generalized Borromean Rings are the collection of links obtained from the Hopf link by (any non-trivial amount of) iterated ramified Bing doubling. 
There is a natural partial order on GBRs where more ramification and more Bing doubling means ``less than'', see figure \ref{GBRs}.
}

An $n$-component link $K$ is called {\em h-trivial}$^+$ if each one of the $n$ links obtained by adding to $K$ a parallel copy of a single component is homotopically trivial. 
The extension from h-trivial to h-trivial$^+$ is of interest in part due to the theorem \cite{FT2} that untwisted Whitehead doubles of h-trivial$^+$ links are topologically slice. This means that  {the strongest possible}  version of the A-B slice problem for such links has a solution. That is, given an h-trivial$^+$ link, there exists a stabilization giving a relatively slice link. It is interesting to relate this to Theorem \ref{homotopy+ theorem} which, starting with a GBR, finds a stabilization yielding an h-trivial$^+$ link. A natural question is whether these two stabilizations can be combined to give a genuine A-B slice solution for GBRs.

As discussed in \cite{FK}, if a link homotopy solution can be sufficiently improved, this could lead to an affirmative resolution to the surgery conjecture for all groups. Our result may be seen as a step in this direction. In terms of Milnor's $\bar\mu$-invariants, \cite{FK} constructs a stabilization giving a link with trivial $\bar\mu$-invariants with non-repeating indices. Theorem \ref{homotopy+ theorem} improves this to a link with trivial $\bar\mu$-invariants with at most two repeating indices.
An interesting question is whether there exists a stabilization giving a link with {\em all} vanishing $\bar\mu$-invariants.

The proof of theorem 1 may be extended to give a link-homotopy$^{+k}$ solution, where $+k$ means that any link obtained by adding a total of $k$ parallel copies of various components is homotopically trivial. However this gain comes at a price: the amount of Bing doubling in GBRs for which our methods give a link-homotopy$^{+k}$ solution grows with $k$.  (The simplest representative link has $2^{2k}+1$ components.) It follows from grope height raising \cite[Proposition 2.7]{FQ} that such a collection of links, for a fixed $k$, is still universal for surgery, see \cite[Proposition 4.1]{FK}.

The novel ingredient in the construction of a link-homotopy solution in \cite{FK} is a geometric use of the group-theoretic $2$-Engel relation $[[y,x],x]$, in conjunction with handle slides. An important algebraic feature underlying this construction is the fact that any $2$-Engel group is nilpotent of a fixed class. 

One natural way to approach the h-trivial$^+$ condition is to use the $3$-Engel relation. However, as discussed in \cite{FK}, $n$-Engel relations for $n>2$ are generally not well understood. Instead, in this paper h-triviality$^+$ is achieved through a systematic application of the $2$-Engel relation. More precisely, we study links modulo a more subtle relation which may be termed ``$(2,2)$-Engel''. 
Roughly, a group element satisfies this relation if it is trivial modulo the $2$-Engel relation in two different ways. h-triviality$^+$ is seen to be a geometric consequence of this relation.

The background material on the A-B slice problem and the $2$-Engel relation is summarized in sections \ref{background section}, \ref{2Engel section}. Theorem \ref{homotopy+ theorem} is proved in section
\ref{proof section}.

\section{A-B slice links and the relative slice problem} \label{background section}

This section gives a brief summary of the relevant background on the A-B slice problem, the reader is referred to \cite{FL}, for a more detailed exposition. We also state the notion of a link-homotopy$^+$ solution to the A-B slice problem, used in theorem~\ref{homotopy+ theorem}.

A  {\em decomposition} of $D^4$ is a pair of compact codimension zero smooth submanifolds with boundary $A,B\subset D^4$, satisfying conditions $(1)$-$(3)$ below. Denote $$\partial^{+} A=\partial A\cap \partial D^4, \; \; \partial^{+} B=\partial B\cap \partial D^4,\; \; \partial A=\partial^{+} A\cup {\partial}^{-}A, \; \; \partial B=\partial^{+} B\cup {\partial}^{-}B.$$ (1) $A\cup B=D^4$,\\ (2) $A\cap B=\partial^{-}A=\partial^{-}B,$ \\
(3) $S^3=\partial^{+}A\cup \partial^{+}B$ is the standard genus $1$
Heegaard decomposition of $S^3$.

The ``attaching curves'' ${\alpha}, {\beta}$ of $A, B$ (the cores of the solid tori $\partial^{+}A$, respectively  $\partial^{+}B$) form the Hopf link in $S^3=\partial D^4$.

Given a $k$-component link $L=(l_1,\ldots,l_k)\subset S^3$, let $L'=(l'_1, \ldots, l'_k)$ be its untwisted parallel copy. 

\begin{definition}  \label{AB slice definition}
A link $L$ is {\em A-B slice} if there exist decompositions $(A_i, B_i)$,
of $D^4$ and self-homeomorphisms ${\phi}_i, {\psi}_i$
of $D^4$, $i=1, \ldots, k$ such that all sets 
${\phi}_1 A_1,  \ldots, {\phi}_k  A_k$,  ${\psi}_1 B_1,\ldots, {\psi}_k
B_k$ are disjoint and satisfy:
${\phi}_i({\partial}^{+}A_i)$  is a tubular neighborhood of $l_i$
and ${\psi}_i({\partial}^{+}B_i)$ is a tubular neighborhood of
$l'_i$, for each $i$.
\end{definition}

The collection of $2k$ manifolds $\{A_i, B_i\}$ are disjointly embedded into $D^4$
by the restrictions ${\phi}_i|_{A_i}$, ${\psi}_i|_{B_i}$.
Since these maps are restrictions of
self-homeomorphisms of $D^4$, the embeddings are {\em standard}, in the sense that 
the complement
$D^4\smallsetminus {\phi}_i(A_i)$ is homeomorphic to $B_i$, and
$D^4\smallsetminus {\psi}_i(B_i)$ is homeomorphic to  $A_i$.
(This  condition that the  embeddings are standard is important: it was shown in \cite{K} that any link with trivial linking numbers is {\em weakly} A-B slice, when this condition is omitted.)
A version of this requirement, in the link-homotopy setting, is stated as condition \ref{standard homotopy definition}.

The $4$-dimensional topological surgery conjecture for free groups was reformulated in \cite{F2, F3}  in terms of the A-B slice problem for the Generalized Borromean rings (GBRs), the collection of links formed from the Borromean rings by iterated ramified Bing doubling. An example is shown in figure \ref{model link}.

The notion of a link-homotopy$^+$ A-B slice link relies on a choice of handle decompositions of the submanifolds $\{ A_i, B_i\}$. We will analyze them in the context of the {\em relative-slice problem}, introduced in \cite{FL}.

Given a decomposition $D^4=A\cup B$, without loss of generality it may be assumed \cite{FL} that each side $A, B$  has a handle decomposition 
(relative to the collar $S^1\times D^2\times I$) with only $1$- and $2$-handles. Denote $A=({\partial}^+ A)\times I\,\cup\, { H}_1\,\cup\, { H}_2$.
As usual, the $1$-handles will be considered as standard $2$-handles $ H_1^*$ removed from the collar, $A=({\partial}^+ A\times I\setminus  { H}^*_1)\cup { H}_2$.
The decompositions constructed in this paper (see section \ref{proof section}) have the property that the $2$-handles $H_2$ of each side do not go through the handles $H_1^*$ of the same side. (See \cite{FK} for a discussion of this terminology.)

Suppose an $n$-component link $L$ is A-B slice, with decompositions $D^4=A_i\cup B_i$, $i=1,\ldots, n$. Denote by $D^4_0$ a smaller $4$-ball obtained by removing from $D^4$ the collars on the attaching regions ${\phi}_i(\partial^+ A_i), {\psi}_i(\partial^+ B_i)$ of all submanifolds $\{ {\phi}_i(A_i)$, ${\psi}_i(B_i)\}$. Let ${\mathcal H}_2$ denote the $2$-handles of all these submanifolds, and ${\mathcal H}_1^*$ the $2$-handles removed from the collars, corresponding to the $1$-handles. Consider ${\mathcal H}_1^*$ as zero-framed $2$-handles attached to $D^4_0$. (See \cite[Section 3.1]{FK} for more details.)

\begin{figure}[ht]
\centering
\includegraphics[height=3.6cm]{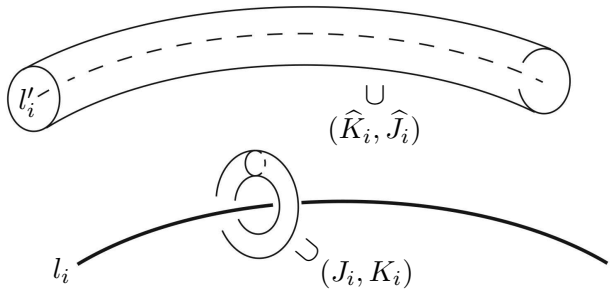}
{\small
    \put(-227,62){$l'_i$}
    \put(-214,-2){$l_i$}
\put(-113,-3){$(J_i, K_i)$}
\put(-110,51){$(\widehat K_i, \widehat J_i)$}}
\put(-96,64){$\cup$}    
\put(-123,5){$\slantedcup$}
 \caption{Stabilization corresponding to an A-B slice link $L=\{l_i\}$: link pairs 
$(J_i, K_i)\subset\, $solid torus neighborhood of a meridian to $l_i$,
$(\widehat K_i, \widehat J_i)\subset \, $solid torus neighborhood of a parallel copy $l'_i$, and 
a diffeomorphism between the solid tori exchanging their meridian and longitude takes $K_i$ to $\widehat K_i$ and $J_i$ to $\widehat J_i$.
}
\label{RSlice figure}
\end{figure}

Consider the following two links $J, K$ in $S^3=\partial D^4_0$.
$J$ denotes the attaching curves of the $2$-handles ${\mathcal H}_2$, 
and $K$ the attaching curves of the $2$-handles ${\mathcal H}_1^*$. (Here ${\mathcal H}_2$ are $2$-handles contained in $D^4_0$, and ${\mathcal H}_1^*$ are attached  to $D^4_0$ with zero framings along $K$.)  We call the pair $(J,K)$ a {\em stabilization} of the original link $L$. The structure of the stabilization links, which is a consequence of the duality between the $1$- and $2$-handles of the two sides of each decomposition, is shown in figure \ref{RSlice figure}. In terms of the notation in the figure, $J=L\cup_i J_i\cup_i \widehat K_i$, $K=\cup_i K_i\cup_i \widehat J_i$.

\begin{definition} \label{rel slice definition}
A link pair $(J, K)$ in $S^3=\partial D^4_0$ is called {\em relatively slice} if the components of $J$ bound disjoint, smoothly  embedded disks in the handlebody
$$ H_K=D^4_0\, \cup\,\text {zero-framed 2-handles attached along}\; K.$$
\end{definition}

If a link $L$ is A-B slice, the associated link pair $(J, K)$ is relatively slice. 
We now turn to the definition of  a link-homotopy$^+$ solution to the  A-B slice problem, referred to in the statement of theorem \ref{homotopy+ theorem}. 

\begin{definition} \label{link homotopy plus ABslice definition} 
A $k$-component link $L$ is {\em link-homotopy$^+$ A-B slice} if there exist decompositions $D^4=A_i\cup B_i$, $i=1,\ldots, k$ and handle decompositions of
the submanifolds $A_i, B_i$ so that the associated relative-slice problem $(J, K)$ has a link-homotopy$^+$ solution. In other words, in the context of definition \ref{rel slice definition} for each component $l$ of $J$ the link $J\cup l'$ bounds disjoint maps of disks $\Delta$ in the handlebody $H_K$. Here $l'$ denotes an untwisted parallel copy of $l$. Moreover, the disks $\Delta$ are subject to condition \ref{standard homotopy definition} below.
\end{definition}

Recall that in the formulation of the A-B slice problem the disjoint embeddings of the manifolds $\{ A_i, B_i\}$ are required to be standard, see the paragraph following definition \ref{AB slice definition}.  We formulate a version of this condition in the link-homotopy setting (see \cite[Section 3.1]{FK} for a more detailed discussion):

\begin{condition} \label{standard homotopy definition} 
Let $S$ be any submanifold in the collection $\{ {\phi}_i(A_i)$, ${\psi}_i(B_i)\}$. 
Then the maps of disks $\Delta$ for the components of $J\cup l'$ corresponding to $S$ do not go through (are disjoint from the co-cores of) the $2$-handles 
attached to $D^4_0$ along the components of $K$ corresponding to the same submanifold $S$. 
\end{condition}

\section{The Milnor group and the $2$-Engel relation} \label{2Engel section}

We start by giving a brief overview of the Milnor group \cite{M}.
Let $\pi$ be a group normally generated by a fixed finite collection of elements $g_1, \ldots,  g_k$. The 
{\em Milnor group} of ${\pi}$, defined with respect to a given normal generating set $\{ g_i\}$,
is defined by
\begin{equation} \label{Milnor group definition}
M{\pi} = {\pi} \, /\, \langle\! \langle \, [g_i, g_i^y] \;\; i=1,\ldots, k,  \; \, y \in G \rangle \! \rangle.
\end{equation}

The Milnor group $ML$ of a link $L$ in $S^3$ is set to be the Milnor group $M{\pi}$ where
${\pi}={\pi}_1(S^3\smallsetminus L)$, defined with respect to meridians to the link components.

Denote by $F_{g_1,\ldots, g_k}$ the free group generated by the $\{g_i\}$, $i=1,\ldots, k$. The Magnus
expansion
\begin{equation}\label{Magnus}
M\co F_{g_1,\ldots, g_k}\longrightarrow {\mathbb Z}[\! [x_1,\ldots, x_k]\! ]
\end{equation}
into the ring of formal power series in non-commuting variables $\{ x_i\}$ is given by $$M(g_i)=1+x_i, \;\,
M(g_i^{-1})=1-x_i+x_i^2-x_i^3\pm\ldots$$ The Magnus expansion induces an injective homomorphism
\begin{equation} \label{MagnusMilnor}
MF_{g_1,\ldots, g_k}\longrightarrow R_{x_1,\ldots,x_k},
\end{equation}
into the quotient $R_{x_1,\ldots,x_k}$ of ${\mathbb Z}[\! [x_1,\ldots, x_k]\!  ]$ by the ideal generated by all monomials
$x_{i_1}\cdots x_{i_k}$ with some index occuring at least twice. Milnor's $\bar\mu$-invariants of a link are defined in terms of coefficients of the Magnus expansion \cite{M2}.

Two links are {\em link-homotopic} if they
are connected by a $1$-parameter family of link maps where different components
stay disjoint for all values of the parameter.
If $L$, $L'$ are link-homotopic then 
their Milnor groups $ML$, $ML'$ are isomorphic. Moreover, a $k$-component link $L$ is
homotopically trivial  (h-trivial) if and only if $ML$ is isomorphic to the free Milnor group $MF_{m_1,\ldots, m_k}$. Equivalently, a link is h-trivial if and only if all its $\bar\mu$-invariants with non-repeating indices are trivial.

This paper concerns a stronger version of this equivalence relation. An $n$-component link $L$ is called {\em h-trivial}$^+$  if each one of the $n$ links obtained by adding to $L$ a parallel copy of a single component is homotopically trivial. A link is h-trivial$^+$  if and only if all its $\bar\mu$-invariants with at most two repeating indices are trivial.

\subsection{2-Engel groups}
Given a group ${\pi}$, consider its lower central series defined by ${\pi}^1={\pi}, {\pi}^n=[{\pi}^{n-1}, {\pi}]$. It is convenient to introduce a concise commutator 
notation $$[g_1, g_2, \ldots, g_n]:= [[\ldots [g_1, g_2],\ldots, g_{n-1}], g_n].$$

This paper concerns geometric applications of the $2$-Engel relation $[[y,x],x]=1$, or equivalently $[x,x^y]=1$. A {\em $2$-Engel group} ${\pi}$ is a group satisfying this relation for all $x,y\in {\pi}$.
Note the difference with the definition of the Milnor group (\ref{Milnor group definition}) where this relation is imposed only on $x$ in a fixed set of normal generators.

The free Milnor group on $n$ generators $MF_n$ is nilpotent of class $n$ \cite{M}. In contrast, the nilpotency class of $2$-Engel groups is independent of the number of generators. This result, building on earlier work of Burnside \cite{Burnside}, is due to Hopkins \cite{Hopkins} (also see \cite{Levi}):

\begin{lemma}  \label{nilpotent} \sl
Any $2$-Engel group is nilpotent of class $\leq 3$.
\end{lemma}

A proof in the context of the Milnor group is given in \cite{FK}. The following more precise statement  will be useful for applications in the next section.

\begin{corollary} \cite[Corollary 2.3]{FK} \label{Engel corollary} \sl
Suppose $\pi$ is a group normally generated by $g_1,\ldots, g_n$. Let $g\in {\pi}^k$ be an element of the $k$-th term of the lower central series, $4\leq k\leq n$. Then $g$ may be represented in the Milnor group $M{\pi}$ as a product of (conjugates of) $k$-fold commutators $C$ of the form $[h_1,\ldots, h_k]$ where two of the elements $h_i$ are equal to each other and to a product of two generators, $h_j = h_m=g_{i_1} g_{i_2}$ for some $j\neq m$, and each other element $h_i$ is one of the generators $g_1, \ldots, g_n$.  
\end{corollary}

We will refer to the commutators of the form $[h_1,\ldots, h_k]$ in the statement of Corollary \ref{Engel corollary} as {\em almost basic commutators}. 
{(This term is meant to avoid confusion with commutators $[g_{i_1},\ldots,g_{i_k}]$ which are usually called {\em basic}.)}

\begin{figure}[ht]
\centering
%\vspace{.5cm} 
\includegraphics[width=12.7cm]{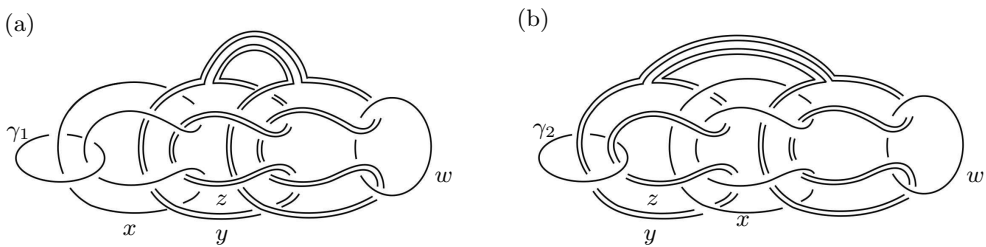}
{\scriptsize
\put(-364,33){${\gamma}_1$}
\put(-320,-4){$x$}
\put(-285,-6){$y$}
\put(-285,7){$z$}
\put(-202,16){$w$}
\put(-165,34){${\gamma}_2$}
\put(-123,-6){$y$}
\put(-1,16){$w$}
\put(-88,-1){$x$}
\put(-122,7){$z$}
\put(-364,73){(a)}
\put(-170,75){(b)}
}
\caption{(a): $\!{\gamma}_1=[x,yz,yz,w]$, $\;$ (b): $\!{\gamma}_2=[yz,x,yz,w]$.}
\label{Elementary links figure}
\end{figure}

Figure \ref{Elementary links figure} shows examples of  links which are a geometric realization (for $k=4$) of almost basic commutators. A central feature of these homotopically {\em essential} links is that a $0$-framed handle slide (in the notation of the figure, of the $z$-curve over the $y$-curve) gives a split link consisting of an unknot and a homotopically {\em trivial} link pictured in figure \ref{Engel links}. These links will be used to construct decompositions of $D^4$ in section \ref{proof section}.

Section \ref{step3} will use the following variation of the statement of Corollary \ref{Engel corollary}.

\begin{proposition} \label{5fold} {\sl 
Let $\pi$ be a group normally generated by $g_1,\ldots, g_n$. Suppose $g\in {\pi}^k$, $5\leq k\leq n$. 
Then $g$ may be represented in the Milnor group $M{\pi}$ as a product of (conjugates of) $k$-fold almost basic commutators $[h_1,\ldots, h_k]$,
where the first element $h_1$ is a generator $g_i$ (and not a product of two generators).  }
\end{proposition}

{\em Proof.} To be specific, consider an element $g\in {\pi}^k$, where $k=5$. (The proof for any $k\geq 5$ is analogous.) Considered as an element of $(MG)^5$, $g$ equals a product of conjugates of $5$-fold basic commutators $[g_{i_1},\ldots,g_{i_5}]$ with distinct indices.
According to Corollary \ref{Engel corollary}, $[g_{i_2}, \ldots,g_{i_5}]$  is a product of conjugates of almost basic commutators $C$. Using commutator identities (\ref{identities}), $[g_{i_1},\ldots,g_{i_5}]$ is then seen to be a product of conjugates of elements of the form $[g_{i_1}, C]$. 
\QEDB

Links corresponding to $5$-fold commutators as in the statement of proposition \ref{5fold} are shown  in figure \ref{Elementary links figure1}. (Note that the component labeled $g_i$ in the figure is not involved in band sums forming the links.)

\section{A motivating example} \label{motivating subsection}
\begin{figure}[ht]
\centering
\includegraphics[width=6.5cm]{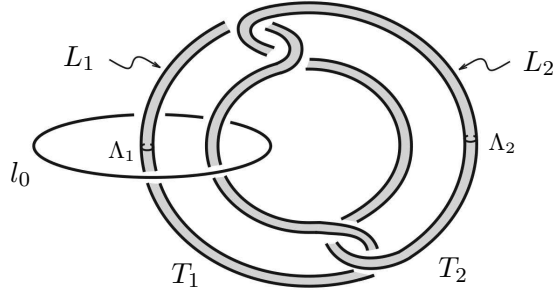}
{\small
\put(-192,42){$l_0$}
\put(-131,4){$T_1$}
\put(-30,6){$T_2$}
\put(-172,86){$L_1$}
\put(2,84){$L_2$}}
\scriptsize{
\put(-155,51){${\Lambda}_1$}
\put(-11,54){${\Lambda}_2$}
}
\caption{The link $L=l_0\cup L_1\cup L_2$ in lemma \ref{motivation lemma}.}
\label{ExampleLink figure}
\end{figure}

Before giving a formal proof of theorem \ref{homotopy+ theorem} we
 illustrate the  idea underlying h-triviality$^+$ in the set-up  in figure \ref{ExampleLink figure}. 
Start with the Borromean rings and let $T_1, T_2$ denote solid torus neighborhoods of two of the components.
Let $L_i\subset T_i$ be two links embedded in these solid tori. Denote by ${\Lambda}_i$ a meridian of $T_i$: a curve in $\partial T_i$ bounding a disk in $T_i$, figure \ref{ExampleLink figure}. 

\begin{lemma} \label{motivation lemma} \sl
Consider the link $L=l_0\cup L_1\cup L_2$, figure \ref{ExampleLink figure}. Suppose that for each $i=1,2$, 

\noindent
(1) $L_i\cup{\Lambda}_i$ is h-trivial, and \nl
(2) $L_i$ is h-trivial$^+$ in $S^3$, where $L_i\subset T_i\subset S^3$ and $T_i\subset S^3$ is the standard inclusion.

Then $L$ is h-trivial$^+$.
\end{lemma}

\begin{figure}[ht]
\centering
\includegraphics[width=6cm]{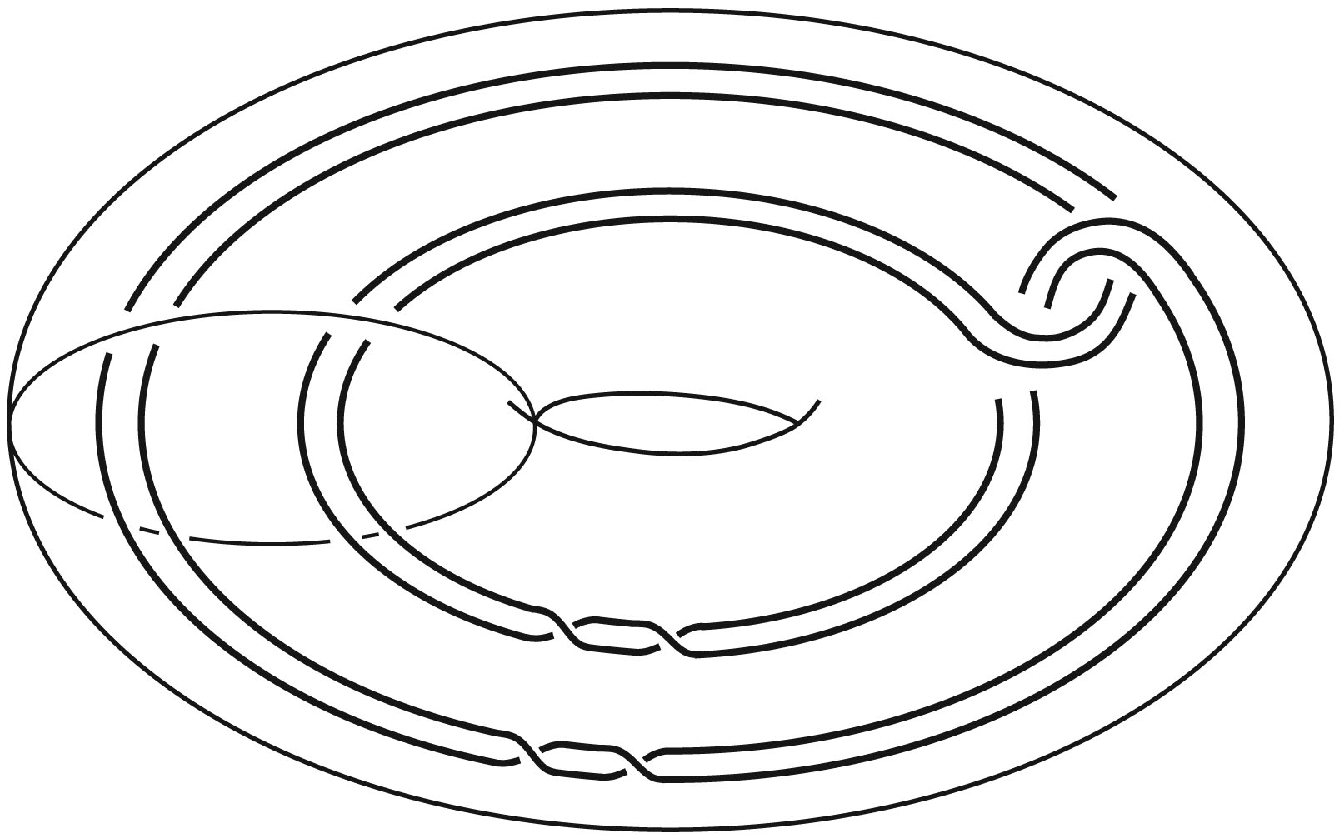}
{\small
\put(-184,49){${\Lambda}_i$}
\put(-39,36){$L_i$}
}
\caption{An example of a link $L_i$ in the solid torus satisfying the assumptions  in lemma \ref{motivation lemma}.}
\label{WhiteheadTorus figure}
\end{figure}

Figure \ref{WhiteheadTorus figure} shows an example of a link $L_i$ in the solid torus satisfying the assumption in lemma \ref{motivation lemma}: the Whitehead double of the core of the solid torus, and a parallel copy. (Note that 
in this example $L_i\cup{\Lambda}_i$ is h-trivial, but $L_i$ is not h-trivial in the solid torus.) Other examples are given by the links in figure \ref{Engel links}, where the three components on the right form the link in the solid torus = complement in $S^3$ of the left-most component.

{\em Proof of lemma} \ref{motivation lemma}. 
There are two separate cases to consider: when a parallel copy of a component of $L_1\cup L_2$ is added, and when a parallel copy of $l_0$ is added.
First consider $L':=L\, \cup\!$ parallel copy $l'_1$ of a component of $L_1$. We start with a geometric argument to show that $L'$ is h-trivial. The steps below are labeled for referencing in follow-up sections.

(1) The link $L$ may be built starting with $L_2\cup{\Lambda}_2$ as follows: Bing double ${\Lambda}_2$, denote one of the resulting components by $l_0$ and insert $L_1$ in a tubular neighborhood of the other component. Pictured this way, $l_0\cup L_1$ is contained in a tubular neighborhood of ${\Lambda}_2$. 
Consider a link null-homotopy of $L_2\cup {\Lambda}_2$ and extend it to $L_2\cup (l_0\cup L_1)$. Self-intersections of ${\Lambda}_2$  during the link-homotopy
are implemented by self-intersections of $l_0$. This gives a link-homotopy of $L$ where the components of $L_1$ have no self-intersections, so the same argument goes through when a parallel copy $l'_1$ is added to $L_1$. The assumption (2) of the lemma completes the proof that $L'$ is h-trivial.
h-triviality of $L\, \cup\!$ parallel copy of a component of $L_2$ is established analogously.

(2) Now we give another, algebraic, proof that $L'$ is h-trivial. This argument will be applicable in the more general setting of theorem  \ref{homotopy+ theorem}.
Abusing the notation, let $l_0, {\Lambda}_i$ refer to based loops. Then 
\begin{equation} \label{first eq}
l_0=[{\Lambda}_1, {\Lambda}_2],
\end{equation}
where
\begin{equation} \label{second eq}
{\Lambda}_1={\Lambda}_2=1\in M{\pi}_1(S^3\smallsetminus (L_1\cup L_2)).
\end{equation}

The statement (\ref{second eq}) follows from the fact that  $(L_1\cup {\Lambda}_1)$ and $(L_2\cup {\Lambda}_2)$ form a split link in $S^3\times 1$ which is homotopically trivial by the assumption (1) of the lemma.

Therefore every monomial (other than $1$) in the Magnus expansion of $l_0$ has two sets of repeated variables: one pair corresponding to a component of $L_1$ and another pair corresponding to a component of $L_2$. This implies that the link remains h-trivial when a parallel copy is added to one component of either $L_1$ or $L_2$.

(3)
Now consider $L\cup l'_0$, where $l'_0$ is a parallel copy of $l_0$.
Note that there exist maps of disks ${\Delta}$ into $D^4$ bounded by $L_1\cup L_2$  and a capped punctured torus $T^c$ bounded 
by $l_0$ in $D^4$ such that all disks and $T^c$ are pairwise disjoint. The body of $T^c$ is an embedded genus $1$ surface bounded by $l_0$ in 
$S^3\smallsetminus (T_1\cup T_2)$, with a symplectic basis of curves isotopic to ${\Lambda}_1, {\Lambda}_2$. Extend $L_1\cup {\Lambda}_1\cup L_2\cup {\Lambda}_2$ by a product in a collar $S^3\times I\subset D^4$, where $S^3\times 0$ is identified with $\partial D^4$.  Then $L_1\cup {\Lambda}_1$ and $ L_2\cup {\Lambda}_2$ form a split link in $S^3\times 1$ which is h-trivial by the assumption (1) of the lemma. The null-homotopies for ${\Lambda}_1, {\Lambda}_2$ give the caps for $T^c$. Contraction/push-off \cite[Section 2.3]{FQ} applied to $T^c$ and its parallel copy give disjoint maps of disks for $l_0, l'_0$ in the complement of the disks $\Delta$ bounded by $L_1\cup L_2$.  This concludes the construction of disjoint disks for all components of $L\cup l'_0$. \QEDB

{\em Remark.} Contraction/push-off in part (3) of the proof, if desired, could be iterated to show that $L\,\cup\!$ (any given number of parallel copies of $l_0$) is a homotopically trivial link.

\section{A lemma in commutator calculus.}

Lemma \ref{motivation lemma} above illustrates the idea that the improvement from a homotopy solution in \cite{FK} to a homotopy$^+$ solution will follow from allowing two ``parallel channels'' or ``two participants in a commutator''  by which an element can die. More precisely, the key features of the link $l_0\cup L_1\cup L_2$ in figure \ref{ExampleLink figure} are the expression (\ref{first eq}) for $l_0$, subject to (\ref{second eq}), and the fact that $L_1\cup L_2$ is h-trivial$^+$.
The link that will come up in the proof of theorem \ref{homotopy+ theorem} is more general than the basic example in lemma \ref{motivation lemma}. 
This section develops the relevant algebraic framework which generalizes Corollary \ref{Engel corollary} using the main features described above.

Note that the representation of $g$ in Corollary \ref{Engel corollary} as a product of conjugates of almost basic commutators holds {\em in the Milnor group}, in general it is not valid in the group ${\pi}$. For the purpose of proving  Theorem \ref{homotopy+ theorem} it is insufficient to work modulo the Milnor relation, a more subtle equivalence relation is needed. 

Let ${\pi}$ be a group normally generated by a fixed set of elements $\{ g_1,\ldots, g_n\}$. Motivated by the defining relations (\ref{Milnor group definition}) of the Milnor group and by the equations (\ref{first eq}), (\ref{second eq}), consider commutators of the form
\begin{equation} \label{22Milnor definition}
\big[ [g_i,g_i^{y_1}]^{z_1}, [g_j,g_j^{y_2}]^{z_2}\big],
\end{equation}

where $1\leq i,j\leq n$, and $y_k, z_k$ are arbitrary elements of $\pi$.
The notation $f\equiv g$ for two elements $f,g\in{\pi}$ will indicate that $f\cdot g^{-1}$ is in the normal subgroup generated by the elements (\ref{22Milnor definition}).
 Technically we will not consider quotients of groups by these relations, rather (\ref{22Milnor definition}) will be used in section \ref{proof section} to construct specific h-trivial$^+$ links.

The following lemma establishes a version of Corollary \ref{Engel corollary} in the setting of the relations (\ref{22Milnor definition}). A useful fact about the lower central series to keep in mind is that $[{\pi}^p,{\pi}^q]\subset {\pi}^{p+q}$.

\begin{lemma} \label{useful lemma1} \sl
Let ${\pi}$ be a group normally generated by $\{ g_1,\ldots, g_n\}$.  Fix $k\geq 4$ and consider a commutator $ [{\alpha},{\beta}]$
where ${\alpha},{\beta}$  are both elements of the $k$th term of the lower central series ${\pi}^k$. Then there exists $W\in{\pi}^{2k}$ such that

(1)   $[{\alpha},{\beta}] \equiv W$, and

(2) $W$ equals in ${\pi}$ a product of conjugates of elements of the form 
$[C, {\beta}]$, $[{\alpha}, C']$, $[C,C']$ where $C, C'$ are almost basic commutators $[h_1,\ldots, h_k]$ of the form introduced in Corollary \ref{Engel corollary}. 
Moreover, for $k\geq 5$ the almost basic commutators $C, C'$ may be assumed to be of the form introduced in Proposition \ref{5fold}.
\end{lemma}

{\em  Proof of lemma \ref{useful lemma1}.} 
Let ${\alpha}'$, ${\beta}'\in {\pi}$ denote  products of conjugates of almost basic commutators representing ${\alpha}, {\beta}$, given by Corollary \ref{Engel corollary}, 
\begin{equation} \label{alphaprime}
{\alpha}\cdot({\alpha}')^{-1}=1 \in M{\pi}, \; \, {\beta}\cdot({\beta}')^{-1}=1 \in M{\pi}.
\end{equation}
Note that in general this does not imply $[{\alpha},{\beta}]\equiv [{\alpha}',{\beta}']$.
Recall the basic commutator identities, cf. \cite[Theorem 5.1]{MKS}:
\begin{equation} \label{identities}
[x,yz]=[x,z]\cdot[x,y]^z, \; \, [xz,y]=[x,y]^z\cdot [z,y].
\end{equation}
It follows from (\ref{alphaprime}) that ${\alpha}\cdot({\alpha}')^{-1}$, ${\beta}\cdot({\beta}')^{-1}$ equal products of conjugates of the defining Milnor relations 
(\ref{Milnor group definition}):
$$ {\alpha}\cdot({\alpha}')^{-1}\, =\, \prod_i \, [g_i,g_i^{x_i}]^{y_i}, \;\,  {\beta}\cdot({\beta}')^{-1}\, =\, \prod_j \, [g_j,g_j^{z_j}]^{w_j},
$$
where $x_i, y_i, z_j, w_j$ are arbitrary elements of $\pi$.
Substituting these expressions for ${\alpha}\cdot ({\alpha}')^{-1}, ({\beta}')^{-1}\cdot {\beta}$ and
using (\ref{identities}), the commutator $[{\alpha}\cdot ({\alpha}')^{-1}, ({\beta}')^{-1}\cdot {\beta}]$ is seen to equal a product of conjugates of elements of the form (\ref{22Milnor definition}).
In other words, 
\begin{equation} \label{eq1}
[{\alpha}\cdot ({\alpha}')^{-1}, ({\beta}')^{-1}\cdot {\beta}] \equiv   1.
 \end{equation}

Again using (\ref{identities}), 
\begin{equation} \label{eq2}
 [{\alpha}\cdot ({\alpha}')^{-1},  ({\beta}')^{-1}\cdot {\beta}]= 
[{\alpha}, {\beta}]^{{\gamma}_1} \cdot 
[({\alpha}')^{-1}, {\beta}] \cdot
[{\alpha}, ({\beta}')^{-1}]^{\gamma_2} \cdot
[({\alpha}')^{-1}, ({\beta}')^{-1}]^{\gamma_3}
\end{equation}

for some ${\gamma}_i\in {\pi}$, determined by (\ref{identities}). Set $W\in {\pi}$ to be (a conjugate of) the inverse of the product of the three right factors in (\ref{eq2}):
\begin{equation} \label{def W}
W:=\big( \big( [({\alpha}')^{-1}, {\beta}] \cdot
[{\alpha}, ({\beta}')^{-1}]^{\gamma_2} \cdot
[({\alpha}')^{-1}, ({\beta}')^{-1}]^{\gamma_3} ]\big)^{-1}\big)^{{\gamma_1}^{-1}}.
\end{equation}

It follows from (\ref{eq1}), (\ref{eq2}) that 
\begin{equation} \label{eq3}
[{\alpha}, {\beta}]\equiv W.
\end{equation}

An application of the basic commutator identities \cite[Theorem 5.1]{MKS},
\begin{equation} \label{identities2}
{[x,y^{-1}]=[y,x]^{y^{-1}},\; [x^{-1},y]=[y,x]^{x^{-1}},}
\end{equation}
allows one to eliminate the negative exponents of ${\alpha}'$, ${\beta}'$ in the definition (\ref{def W}) of $W$, thereby expressing $W$ in the form required in the statement (2) of the lemma.
\QEDB

The proof of theorem \ref{homotopy+ theorem} in section \ref{proof section} will require the following extension of lemma \ref{useful lemma1}.

\begin{lemma} \label{extension} 
Let ${\pi}$ be a group normally generated by $\{ g_1,\ldots, g_n\}$.  Fix $k\geq 4$, and suppose $g\in{\pi}$ is of the form
\begin{equation} \label{lemma eq}
g=\big[ [{\alpha}_1,{\beta}_1],[{\alpha}_2,{\beta}_2]\big],
\end{equation} 
where each element ${\alpha}_i,{\beta}_i$  is in ${\pi}^k$. 
Then there exist $W_1, W_2\in {\pi}^{2k}$, satisfying the condition (2) in Lemma \ref{useful lemma1}, such that $g\equiv [ W_1,W_2]$. 
\end{lemma}

{\em Proof.}
An application of lemma \ref{useful lemma1} to $[{\alpha}_i, {\beta}_i]$ gives  $W_1,W_2\in{\pi}^{2k}$ such that $[{\alpha}_i, {\beta}_i]\equiv W_i$, $i=1,2$. 

Observe that given $x,y,z\in {\pi}$, $x\equiv y$ implies $[x,z]\equiv [y,z]$. Indeed,  $x\equiv y$ means that  $x\cdot y^{-1}$ is in the normal subgroup generated by the elements (\ref{22Milnor definition}). Denote this normal subgroup by $N$. Then $x^z\cdot (y^z)^{-1}=(xy^{-1})^z\in N$, so $x^z\equiv y^z$. Therefore $[x,z]\cdot [y,z]^{-1}=x\cdot (x^{-1})^z\cdot (y\cdot (y^{-1})^z)^{-1}=x\cdot (x^{-1})^z\cdot y^z\cdot y^{-1}$ is also in $N$.

It follows that
\begin{equation} \label{g}
g=\big[ [{\alpha}_1,{\beta}_1],[{\alpha}_2,{\beta}_2]\big]\equiv \big[ W_1,[{\alpha}_2,{\beta}_2] \big] \equiv \big[ W_1,W_2 \big].
\end{equation}
\QEDB

\section{Proof of theorem \ref{homotopy+ theorem}} \label{proof section}

\begin{figure}[ht]
\centering
\includegraphics[height=3.2cm]{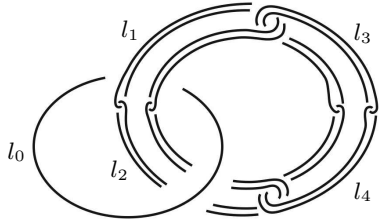}
{\scriptsize
\put(-146,27){$l_0$}
\put(-103,73){$l_1$}
\put(-107,20){$l_2$}
\put(-15,72){$l_3$}
\put(-15,12){$l_4$}
}
 \caption{A link in the collection of GBRs. 
}
\label{model link}
\end{figure}

As discussed in the introduction (also see \cite[Proposition 4.1]{FK}) any coinitial subset of the  Generalized Borromean Rings forms a collection of links 
universal for surgery.
A homotopy A-B slice solution in \cite{FK} applies to links obtained from the Hopf link by keeping one of its components $l_0$ intact and Bing doubling the other components at least twice, see figure \ref{model link} for an example of such a link. A homotopy$^+$ solution constructed in sections \ref{step1} - \ref{step3} applies to
a collection of higher Bing-doubled links (still universal for surgery). 

The decompositions $D^4=A_i\cup B_i$ for all link components other than the fixed component $l_0$ are set to be the trivial decomposition, $A_i=2$-handle and $B_i=$ collar on the attaching curve ${\beta}_i$. Similarly to \cite{FK}, the entire complexity of the construction is in the decomposition $D^4=A_0\cup B_0$ for the component $l_0$. $A_0$ will be defined to be the collar $l_0\times D^2\times I$ with a single $2$-handle attached to the core of the solid torus $l_0\times D^2\times\{ 1\}$, and many  $1$-handles governed by the algebraic outcome of lemma \ref{useful lemma1}, as explained below. In terms of figure \ref{RSlice figure} (where the index $i$ is understood to equal $0$) the $2$-handle is attached to the curve labeled $l_i$, $J_i$ is empty since there are no other $2$-handles, and $K_i$ is the Kirby diagram representation of the $1$-handles.
Correspondingly, the other side $B_0$ of the decomposition has no $1$-handles ($\widehat J_i$ in the figure is empty) and the attaching curves  of its $2$-handles form the link
$\widehat K_i$.

\begin{figure}[ht]
\centering
%\vspace{.5cm} 
\includegraphics[width=12.5cm]{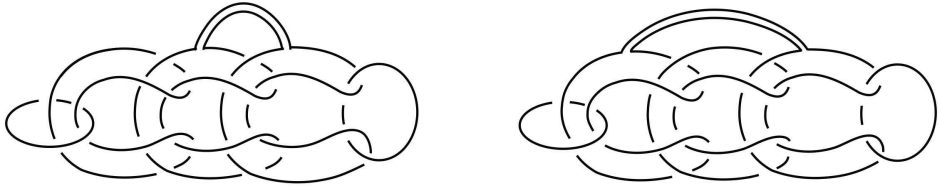}
\caption{Links in figure \ref{Elementary links figure}, with one of the parallel components $y,z$ removed.}
\label{Engel links}
\end{figure}

The link-homotopy solution in \cite{FK} uses a geometric implementation of Corollary \ref{Engel corollary}, where each almost basic commutator of the form $[h_1,\ldots, h_k]$ is realized by a standard link illustrated in figure \ref{Elementary links figure}. 
More precisely, building blocks in the construction of the link $K_0$, describing the $1$-handles of $A_0$, are shown in figure \ref{Engel links}. These are {\em h-trivial}  counterparts of the links in figure \ref{Elementary links figure} where one of the two parallel curves labeled $y,z$ is removed. 

A key point, using the terminology of definition \ref{rel slice definition},  is that $K_0$ is the attaching link for $0$-framed $2$-handles attached to $D^4_0$, and parallel copies of each component bound disjoint copies of the core of the attached $2$-handle. The links in figure \ref{Elementary links figure}  then may be recovered from links in figure \ref{Engel links} by adding the relevant parallel copy.

The homotopy solution from \cite[Section 4]{FK} is schematically shown in Figure \ref{ABslicing} (a). 
A modification of this construction, needed for the proof of theorem \ref{homotopy+ theorem}, is introduced in section \ref{step1}. The theorem is then proved in three steps. The proof of h-triviality$^+$ with respect to the components of $L$, band-summed with parallel copies of $K_0$, relies on lemma \ref{useful lemma1}. h-triviality$^+$ with respect to the components of $\widehat K_0$ uses lemma \ref{extension}. Finally, h-triviality$^+$ with respect to $l_0$ is analyzed in section \ref{step3}.

\subsection{Construction of the stabilization.} \label{step1}
\begin{figure}[ht]
\centering
\includegraphics[height=2.8cm]{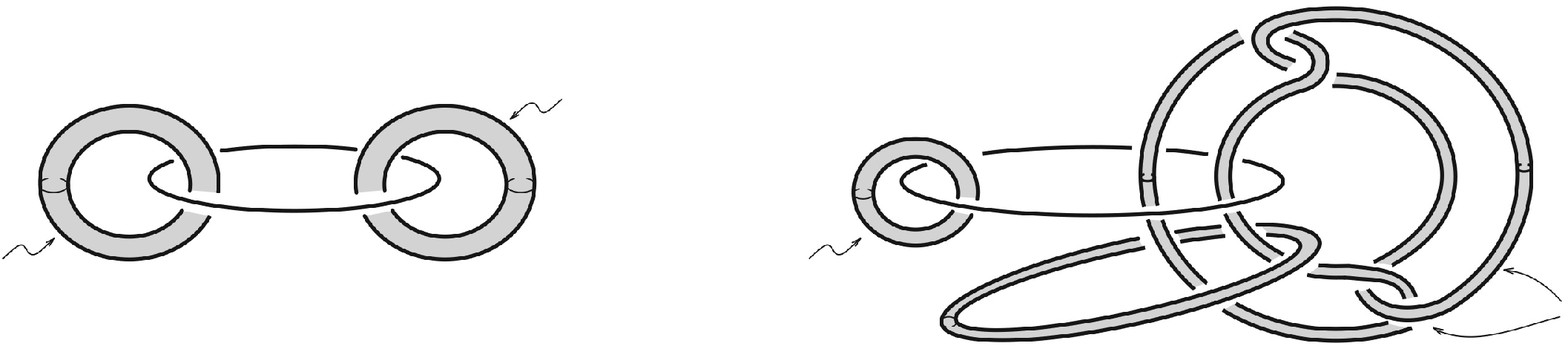}
\small{
\put(-358,65){(a)}
\put(-180,65){(b)}}
{\scriptsize
\put(-356,35){${\alpha}_0$}
\put(-230,34){${\alpha}$}
\put(-295,21){$l_0$}
\put(-357,14){$K_0$}
\put(-250,61){$(l_1,\ldots,l_4)$}
\put(-180,12){$K_0$}
\put(-175,33){$M$}
\put(-154,-2){$\widehat K_0$}
\put(-117,62){$[{\alpha}_1, {\beta}_1]$}
\put(-22,71){$[{\alpha}_2,{\beta}_2]$}
\put(-115,35){$l_0$}
\put(-1,5){$L$}
}
 \caption{\nl
(a) A schematic representation of the link-homotopy A-B slicing in [7], where $(l_0,\ldots,l_4)$ is the GBR is figure \ref{model link}. (The link $\widehat K_0$, and also bands connecting $l_2,\ldots, l_4$ with $K_0$ are not shown.)
\nl 
(b) The modified stabilization $K_0$, constructed in section \ref{step1}, for the GBR $(l_0, L)$. Here $L$ is obtained by Bing doubling twice the components $l_1,\ldots, l_4$ of the link in figure \ref{model link}.}
\label{ABslicing}
\end{figure}

Let {{${\alpha}_1, {\beta}_1,{\alpha}_2, {\beta}_2$ denote meridians to the components $l_1,\ldots, l_4$ of the link in figure \ref{model link}.}}
Consider the link $(l_0, L)$ in figure \ref{ABslicing} (b),  obtained by  Bing doubling twice the components $l_1,\ldots, l_4$. 
{Then ${\alpha}_i, {\beta}_i$ represent $4$-fold commutators in the meridians of $L$. The commutators $[{\alpha}_1, {\beta}_1]$, $[{\alpha}_2, {\beta}_2]$ are meridians of the solid tori as indicated in figure \ref{ABslicing} (b).}
The argument applies to GBRs that are more Bing-doubled and ramified; to be concrete we focus here on this simplest representative link.
This section constructs:
\begin{itemize}
\item{a stabilization of $l_0$: a link $K_0$ in a solid torus linking $l_0$, and }
\item{a band-sum $L^{\sharp}$ of $L$ with ($K_0\, \cup\!$ parallel copies).}
\end{itemize}
The following properties will be established:
\begin{itemize}
\item{The link $K_0$ is h-trivial$^+$ in the solid torus,}
\item{$(l_0\cup L^{\sharp})$  is h-trivial$^+$ with respect to adding a parallel copy of a component of $L^{\sharp}$.}
\end{itemize}

{The remaining step, homotopy triviality$^+$ with respect to a parallel copy of $l_0$ will be established in section \ref{step3}.}

To begin the construction, consider ${\pi}:={\pi}_1(S^3\smallsetminus L)$. $L$ is the unlink, ${\pi}$ is the free group, and $l_0$ is an element of the form 
$\big[ [{\alpha}_1,{\beta}_1],[{\alpha}_2,{\beta}_2]\big]$ as in equation (\ref{g}). In the concrete example discussed above, each ${\alpha}_i, {\beta}_i$ is a $4$-fold commutator. (More generally, the argument applies to higher Bing doubles where ${\alpha}_i, {\beta}_i$ are $k$-fold commutators, $k\geq 4$. To obtain the argument in full generality, one replaces ``$4$'' (respectively ``$3$'') with ``$k$'' (respectively ``$k-1$'') in each reference to the length of commutators or the number of link components in this section.)
According to the statement (2) of lemma \ref{useful lemma1}, each element $W_1, W_2$ in (\ref{g}) is a product of conjugates of the elements of the form 
\begin{equation} \label{terms} 
[C, {\beta}_i], \; [{\alpha}_i, C'], \; [C,C'],
\end{equation}
where $C, C'$ are almost basic commutators as in Corollary \ref{Engel corollary}. 
The stabilizing link $K_0$ will be constructed as a geometric analogue of the commutator $[W_1, W_2]$ in (\ref{g}).
Specifically, consider the Bing double of a meridian to $l_0$, and thicken the two resulting components to solid tori $T_1, T_2$. The link in the solid torus $T_i$ will be constructed to correspond to the commutator $W_i$, $i=1,2$. In the concrete example under discussion, the link $L$ is symmetric and the links constructed in $T_1, T_2$ are going to be identical. (In general, these links depend on the elements $W_1, W_2$  in section \ref{extension}.)

Consider several copies of the core of each solid torus $T_1, T_2$, one for each factor in the statement (2) of lemma \ref{useful lemma1}.
For each element of the form (\ref{terms}) in the expression for $W_i$, 
take a Bing double of the corresponding copy of the core of $T_i$ and thicken the two resulting components to solid tori. Next we define links $K', K''$, geometrically representing the  given elements (\ref{terms}) and insert them into these solid tori. 
For each element of the form $[C,{\beta}_i]$ consider a link of the type shown in figure \ref{Engel links}, corresponding to the almost basic commutator $C$. More precisely, $K'$ consists of three components on the right in a link in figure \ref{Engel links}, considered in the solid torus complement of the leftmost component. Let $K''$ be the iterated Bing double of the core of the solid torus corresponding to the $4$-fold commutator ${\beta}_i$. The analogous pair of links ($K', K''$) is created for each factor of the form (\ref{terms}) in the expression for $W_i$, completing the construction of $K_0$.

The construction in the preceding paragraph is a generalization of that in \cite[Section 4]{FK}, in particular see figure 4.3 in that reference. The stabilization in \cite{FK} 
was defined in terms of almost basic commutators $C$, while here we have links corresponding to commutators (\ref{terms}). 

{Now we construct the band sum $L^{\sharp}$, as promised at the beginning of section \ref{step1}.
For each constituent link of $K_0$ of the type in figure \ref{Engel links} add a parallel copy to one of the components to recreate a link as in figure \ref{Elementary links figure}. In the relative slice setting these parallel copies bound disjoint disks in the zero-framed $2$-handles attached to $D^4_0$ along $K_0$. To find a homotopy$^+$ solution to the relative-slice problem, $L$ will be band-summed with the components of $K_0$ and their parallel copies. 
In the homotopy solution in \cite{FK} the choice of bands was immaterial. This was due to the fact that all commutators in the construction were of maximal length, so conjugation did not affect calculations in the Milnor group.  The only relevant constraint for a homotopy solution was $l_0={\alpha}\cdot ({\alpha}_0)^{-1}=1\in M{\pi}_1(S^3\smallsetminus L)$,  in the notation of figure \ref{ABslicing} (a).  (The equality $l_0={\alpha}\cdot ({\alpha}_0)^{-1}$ was established using additivity of $\bar\mu$-invariants \cite{Cochran, K2}, or by directly reading off the element represented by $l_0$ in the Milnor group \cite[Proof of theorem 1]{FK}.)   The homotopy$^+$ problem is more sensitive to the choice of bands. 

Suppose the bands in the definition of  $L^{\sharp}$ could be chosen so that the meridian $M$ of the solid torus in figure \ref{ABslicing} (b) {\em precisely} matched the element $[W_1,W_2]$ in (\ref{g}). Then in ${\pi}$ one would have $l_0=\big[ [{\alpha}_1,{\beta}_1],[{\alpha}_2,{\beta}_2]\big]\cdot [W_1,W_2]^{-1}\equiv 1$, see (\ref{g}). In this case the h-triviality$^+$ with respect to parallel copies of $L^{\sharp}$ is proved exactly the same way as in part (2) of the proof of lemma \ref{motivation lemma}. 

In fact, to establish h-triviality$^+$  with respect to parallel copies of $L^{\sharp}$ it suffices to choose bands so that $M$ suitably  approximates $W$.
Recall that each $W_1, W_2$ is a product of conjugates of the elements (\ref{terms}). The Milnor group $M{\pi}=M{\pi}_1(S^3\smallsetminus L)$ is nilpotent of class equal to the number of components of $L$. For any choice of bands, each commutator  (\ref{terms}) in the expression for $W_i$ is of maximal length in the power series of non-repeating monomials. Considering the Magnus expansion (\ref{Magnus}) of the free group $\pi$, note that only the homology class of conjugating elements is relevant. Indeed, suppose elements $W_1', W_2'$ are created by some conjugating elements that agree homologically with those defining $W_i$ in Lemma \ref{useful lemma1}. Then each term (other than $1$) in the Magnus expansion of $[W_1,W_2]\cdot [W_1',W_2']^{-1}$  contains either three copies of a variable, or two pairs of different repeated variables. In either case adding a parallel copy preserves the condition of being h-trivial.
Finally, choose arcs connecting $L$ with the relevant components of $K_0$ and its parallel copies, homologically matching the conjugating elements. Perform the band sums along these arcs; such operations do not interfere with each other since only homological information is relevant. This establishes h-triviality$^+$ with respect to components of $L^{\sharp}$.

It remains to show that the link $K_0$ is h-trivial$^+$ in the solid torus (or equivalently that $\widehat K_0$ is h-trivial$^+$ in the solid torus).
Note that each link in the construction above, corresponding to an almost basic commutator $C$, is h-trivial in the solid torus. These links consist of three components on the right in a link in figure \ref{Engel links}, considered in the solid torus complement of the leftmost component. The construction of $K_0$ was based on the commutator $[W_1,W_2]$ in (\ref{g}), and both $W_1, W_2$ are products of conjugates of the elements of the form (\ref{terms}). 
{$K_0$ was defined as the union of links in the two solid tori $T_1, T_2$ (whose cores are the Bing double of a meridian to $l_0$). The link in $T_i$ for each $i=1,2$, is h-trivial  in the solid torus, therefore 
$K_0$ is h-trivial$^+$ in the solid torus.

\subsection{h-triviality$^+$ with respect to ${\mathbf l_0}$} \label{step3}
Finally, consider the link $l_0\cup l'_0\cup L^{\sharp}$ in the notation of section \ref{step1}. The goal is to show that this link is h-trivial. To achieve this, the stabilization will require a slight modification, based on the observation in remark \ref{5fold}.
This argument is independent of section \ref{step1}; it can be carried out in the setting of \cite{FK}.

Consider the links in figure \ref{Elementary links figure1}. They are obtained from the elementary Engel links in figure \ref{Elementary links figure} by adding an additional component labeled $g_i$ in the figure. These links are geometric realizations of $5$-fold almost basic commutators $[h_1,\ldots, h_5]$ in the statement of corollary \ref{Engel corollary}, where the first element $h_1$ is a generator $g_i$, rather than a product of two generators. Geometrically this is manifested in the fact that the band sums used to construct the links do not involve the component labeled $g_i$. As pointed out in remark \ref{5fold}, any element $g\in {\pi}^k$, $k\geq 5$, may be represented in the Milnor group $M{\pi}$ as a product of conjugates of commutators of this type. 

\begin{figure}[ht]
\centering
%\vspace{.5cm} 
\includegraphics[width=12.7cm]{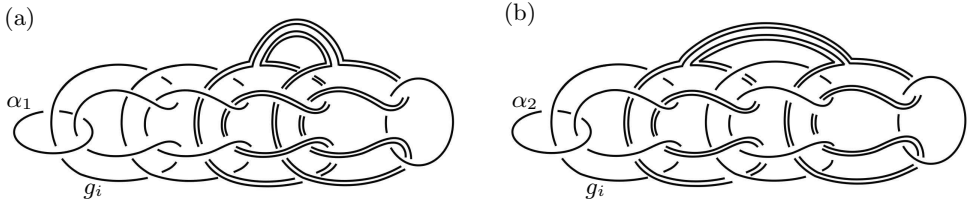}
\scriptsize
\put(-335,-1){$g_i$}
\put(-145,-1){$g_i$}
\put(-364,32){${\alpha}_1$}
\put(-173,32){${\alpha}_2$}
\put(-364,63){(a)}
\put(-175,65){(b)}
\caption{Links in figure \ref{Elementary links figure} with an extra component.
(a): $ {\alpha}_1=[g_i,{\gamma}_1]$, $\;$ (b): ${\alpha}_2=[g_i, {\gamma}_2]$.}
\label{Elementary links figure1}
\end{figure}

To illustrate the use of the links in figure \ref{Elementary links figure1}, first consider the Borromean rings with a parallel copy $l'_0$ of one of its components $l_0$, figure \ref{calc figure}. Denoting by $m_1,m_2$ the meridians suitably connected to a basepoint and similarly regarding $l_0$ as a based loop, note that the expression 
$l_0=[m_1,m_2]$ holds regardless of whether $l_0$ is considered as an element of the Milnor group $M{\pi}_1(S^3\smallsetminus (l_1\cup l_2))$ or as an element
of $M{\pi}_1(S^3\smallsetminus (l'_0\cup l_1\cup l_2))$. This expression can be read off from the torus bounded by $l_0$ in the complement of $l_1, l_2$, figure \ref{calc figure}.

\begin{figure}[ht]
\centering
\includegraphics[width=4.5cm]{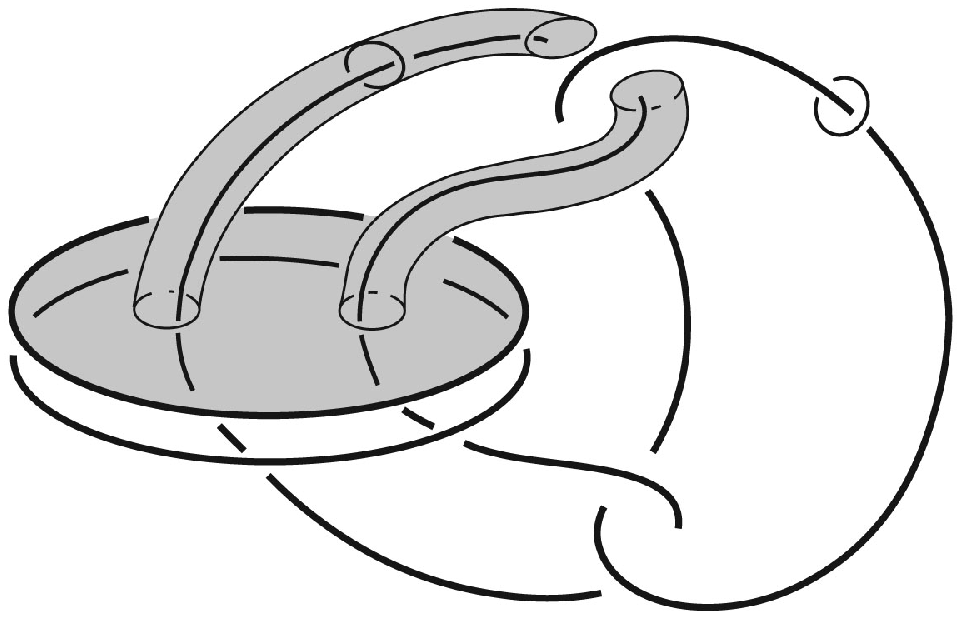}
{\small
\put(-14,71){$m_2$}
\put(-98,78){$m_1$}
\put(-90,0){$l_1$}
\put(-14,0){$l_2$}
\put(-133,46){$l_0$}
\put(-133,24){$l'_0$}
}
\caption{}
\label{calc figure}
\end{figure}

Similarly, the expressions read off by the components ${\alpha}_i$ of the  links in figure \ref{Elementary links figure1} are unchanged when a parallel copy of ${\alpha}_i$ is present. This calculation may be given in terms of the Milnor group, or it can be read off from a grope bounded by $l_0$ in the complement of the link. 
The links in question may be represented as a ``composition'', in the sense of \cite[Theorem 2.3]{FL},  where the components labeled ${\alpha}$ and $g_i$ of a link in figure \ref{Elementary links figure1}  are identified with the components $l_0$, $l_1$ of the Borromean rings in figure \ref{calc figure}, and the rest of the link (denote it by $Q$) is inserted in a solid torus neighborhood $T$ of $l_2$. (This may also be understood as an instance of the more general satellite construction.)  The following argument applies to any link $Q$ in $T$.
As in figure \ref{calc figure}, 
\begin{equation} \label{expression}
l_0=[m_1,m_2].
\end{equation}
Consider the meridian $m_2$ as an element  in the Milnor group of the complement of $Q$ in the solid torus $T$, $M{\pi}_1(T\smallsetminus Q)$. The generators of this Milnor group are meridians to $Q$ and a longitude of the solid torus. Considered as part of figure \ref{calc figure}, this longitude represents the commutator $[m_1, m_0 m'_0]$.
Substituting this into (\ref{expression}), observe that the only way that the meridian $m'_0$ appears in the expression for $l_0$ is as part of the commutator
$$ [m_1,\, \ldots \cdot [m_1,m_0m'_0]\cdot \ldots].$$ 
Applying the commutator identities (\ref{identities}) and the Milnor relation (\ref{Milnor group definition}), it follows that omitting $[m_1,m_0m'_0]$ from this expression does not change 
the element in the Milnor group, establishing the desired claim.

 Now consider a GBR of the form $l_0\cup \overline L$, where $l_0$ is a $5$-fold commutator in the complement of $\overline L$.
(For example in figure \ref{model link} Bing double any one of the components $l_2,\ldots, l_4$.) 
Using remark \ref{5fold}, the construction in section \ref{step1} goes through without any changes, except that  the elements ${\alpha}_i, {\beta}_i$, $i=1,2$ in the link $l\cup L$ in figure \ref {ABslicing} (b) are $5$-fold (rather than $4$-fold) commutators.
{The elementary Engel links used to build $K_0$ are now of the form in figure \ref{Elementary links figure1} (with a parallel copy of the band-summed curves removed), rather than figure \ref{Engel links}.

{Given such GBR $l_0\cup L$ and stabilization $K_0$,  $l_0$ then may be assumed to be in the subgroup (isomorphic to $M{\pi}_1(S^3\smallsetminus (L\cup K_0))$) of $M{\pi}_1(S^3\smallsetminus (l'_0\cup L\cup K_0))$, generated by meridians to $L\cup K_0$. 
As in the proof in section \ref{step1}, only the homological information about
bands connecting $L$ and $K_0$ and forming $L^{\sharp}$ is relevant. Choose arcs connecting $L$ with $K_0$ so that the corresponding conjugating elements are in the subgroup $M{\pi}_1(S^3\smallsetminus (L\cup K_0))$ of $M{\pi}_1(S^3\smallsetminus (l'_0\cup L\cup K_0))$. }
Now all calculations in the preceding sections, establishing that $l_0=1\in M{\pi}_1(S^3\smallsetminus L^{\sharp})$, do not involve the meridian to $l'_0$, so 
 \begin{equation} \label{Milnor calculation}
 l_0=1\in M{\pi}_1(S^3\smallsetminus (l'_0\cup L^{\sharp})).
 \end{equation} 
{The link $l'_0\cup L^{\sharp}$ is h-trivial, therefore its Milnor group is the free Milnor group, and (\ref{Milnor calculation}) is equivalent to h-triviality of the link $l_0\cup l'_0\cup L^{\sharp}$.
 This concludes the proof of theorem \ref{homotopy+ theorem}. \QEDB

{\bf Acknowledgments}.  
VK was partially supported by NSF grant DMS-1309178.
He also would like to thank the IHES for hospitality and support (NSF grant 1002477).

\end{document}